	\documentclass[12pt, a4paper]{article}
	
 	\usepackage{latexsym}
	\usepackage{amsmath}	
	\usepackage{epsfig}	
	\usepackage{times}	
	\usepackage{amssymb}
	
 	\newtheorem{thr}{Theorem}

	\newtheorem{cor}{Corollary}
	\newtheorem{lem}{Lemma}
	
 	\newcommand{\F}{\mathbb{F}}

\begin{document}
 	\centerline{\Large{\bf Divisibility of Trinomials by Irreducible}}
	\centerline{}
	\centerline{\Large{\bf Polynomials over $\F_2$}}
 	\centerline{}
	\centerline{\textsuperscript{a}Ryul Kim, \textsuperscript{b}Wolfram Koepf}
 	\centerline{}
 	\small \centerline{\textsuperscript{a}Faculty of Mathematics and Mechanics, \textbf{Kim Il Sung} University, D.P.R Korea}
	\small \centerline{\textsuperscript{b}Department of Mathematics, University of Kassel, Kassel, Germany}
	\small \centerline{\textsuperscript{a}ryul\_kim@yahoo.com, \quad \textsuperscript{b}koepf@mathematik.uni-kassel.de}
	%\small \centerline{\textsuperscript{b}koepf@mathematik.uni-kassel.de}
	\centerline{}
	\centerline{}
	\begin{abstract}
	Irreducible trinomials of given degree $n$ over $\F_2$ do not always exist and in the cases that there is no irreducible 
	trinomial of degree $n$ it may be effective to use trinomials with an irreducible factor of degree $n$. In this paper we 
	consider some conditions under which irreducible polynomials divide trinomials over $\F_2$. A condition for divisibility 
	of self-reciprocal trinomials by irreducible polynomials over $\F_2$ is established. And we extend Welch's criterion for 
	testing if an irreducible polynomial divides trinomials $x^m+x^s+1$ to the trinomials $x^{am}+x^{bs}+1$.
	\end{abstract}
	%
	%
	%%%%%%%%%%%%%% Section 1 Introduction %%%%%%%%%%%%%%%%
	%
	%
	\section{Introduction}
	Irreducible and primitive trinomials over finite fields are of interest both in theory and practice. 
	We restrict our attention to polynomials over a binary field $\F_2$. Sparse polynomials such as trinomials are 
	commonly used to perform arithmetic in extension fields of finite fields since they provide a fast modular reduction 
	but unfortunately irreducible trinomials of given degree $n$ over $\F_2$ do not always exist. 
	Swan's theorem \cite{swa} rules out $n \equiv 0 \pmod 8$ and also most $n \equiv \pm 3 \pmod 8$. 
	In the cases that there is no irreducible trinomial of given degree $n$, one can always use irreducible polynomials 
	with more than three nonzero terms like pentanomials. But it may be more effective to use (reducible) trinomials with 
	irreducible or primitive factors of degree $n$. In 1994, Tromp, Zhang and Zhao \cite{tro} asked the following question: 
	given an integer $n$, do there exist integers $m, k$ such that
 	\begin{equation*}
	G = \textnormal{gcd} \big( x^m+x^k+1, x^{2^n-1}-1 \big)
	\end{equation*}
	is a primitive polynomial of degree $n$ over $\F_2$? They verified that the answer is yes for $n$ up to 171 and 
	conjectured that the answer is always yes. Blake, Gao and Lambert \cite{bla} confirmed the conjecture for $n \leq 500$ 
	and they also relaxed the condition slightly and asked: do there exist integers $m, k$ such that $G$ has a primitive factor 
	of degree $n$? Motivated by \cite{bla}, Brent and Zimmermann \cite{bre} defined an almost primitive (irreducible) 
	trinomial which is the trinomial with a primitive (irreducible) factor of given degree $n$ and they proposed the algorithms 
	for finding almost primitive (irreducible) trinomials. Doche \cite{doc} called these trinomials (almost irreducible trinomials) 
	as redundant trinomials and gave a precise comparison of running times between redundant trinomials and irreducible 
	pentanomials over finite fields of characteristic 2. In \cite{gol} it was given a positive answer to the latter question and 
	the authors developed the theory of irreducible polynomials which do, or do not, divide trinomials over $\F_2$. They 
	considered some families of polynomials with prime order $p>3$ that do not divide trinomials. 
	To know which irreducible polynomials divide trinomials over $\F_2$ is of interest in many applications such as 
	generation of pseudo-random sequences. In this paper we consider some conditions under which a given irreducible 
	polynomial divides trinomials over $\F_2$. We prove a condition for a given irreducible polynomial to divide 
	self-reciprocal trinomials. 

	Welch’s criterion is a clever one for testing if an irreducible polynomial divides trinomials over $\F_2$. 
	We give a refinement of a necessary condition for divisibility of trinomials $x^{am}+x^{bs}+1$ by a given 
	irreducible polynomial (\cite{che}) and extend Welch's criterion to this type of trinomials.

	%
	%
	%%%%%% Section 2 Divisibility of self-reciprocal trinomials by irreducible polynomials %%%%%%%%%%%%
	%
	%
	\section{Divisibility of self-reciprocal trinomials by irreducible polynomials}
	In this section we consider divisibility of self-reciprocal trinomials by given irreducible polynomials. 
	Let $q$ be a prime power. For a polynomial $f(x)$ of degree $n$ over finite field $\F_q$ a $\it{reciprocal}$ of $f(x)$ is 
	the polynomial $f^*(x)$ of degree $n$ over $\F_q$ given by $f^*(x)=x^n f(1/x)$ and a polynomial $f(x)$ is called 
	$\it{self}$-$\it{reciprocal}$ if $f^*(x)=f(x)$. Numerous results are known concerning self-reciprocal irreducible 
	polynomials over finite fields. In \cite{yuc}, it was studied in detail the order of self-reciprocal irreducible polynomials 
	over finite fields. Let $f \in \F_q[x]$ be a nonzero polynomial with $f(0) \neq 0$. The least positive integer $e$ for 
	which $f$ divides $x^e-1$ is called the $\it{order}$ of $f$ and denoted by $\textnormal{ord}(f)$ [6]. If $f$ is an 
	irreducible polynomial of degree $n$ over $\F_q$ and with $f(0) \neq 0$ then $\textnormal{ord}(f)$ is equal to the order 
	of any root of $f$ in the multiplicative group $\F_{q^n}^*$ and divides $q^n-1$. Below we assume all polynomials 
	to be over $\F_2$. In this case the order of an irreducible polynomial is always odd integer.
	
	In \cite{gol}, it was proved that for prime $p>3$, if there exists a self-reciprocal irreducible polynomial of order $p$ 
	then all irreducible polynomials of the same order do not divide trinomials. In particular, every self-reciprocal 
	irreducible polynomial of prime order $>3$ does not divide trinomials. In fact we can easily see that a self-reciprocal 
	irreducible polynomial $f$ divides trinomials in $\F_2[x]$ if and only if $\textnormal{ord}(f)$ is a multiple of 3. 
	(See exercise 3.93 in \cite{lid}).
	
	Now consider self-reciprocal trinomials. Self-reciprocal irreducible trinomials over $\F_2$ are only of the form 
	$f=x^{2\cdot3^k}+x^{3^k}+1$ which has order $3^{k+1}$. Then which irreducible polynomial divides 
	self-reciprocal trinomials? As above mentioned, the order of a self-reciprocal irreducible polynomial which divides 
	self-reciprocal trinomial is a multiple of 3. Furthermore, we can say a similar thing about the general irreducible 
	polynomials which divide self-reciprocal trinomials. For this we need an auxiliary result.

	%-------------- Lemma 1 --------------
	%
	\begin{lem}
	If an irreducible polynomial $f$ of order $e$ divides a self-reciprocal trinomial $x^{2m}+x^m+1$, then there exists a 
	unique self-reciprocal trinomial of degree $<e$ which is divided by $f$.
	\end{lem}
	\textit{Proof}. Let $\alpha$ be any root of $f$ in a certain extension of $\F_2$ then $\alpha^{2m}+\alpha^m+1$. Write 
 	\begin{equation*}
	m=e \cdot q+r, ~ 0<r<e
	\end{equation*}
	Note that $r$ is not null. If $2r<e$, then $x^{2r}+x^r+1$ is a desired trinomial. Suppose $2r>e$. ($2r \neq e$ because 
	if $2r=e$, then $0=\alpha^{2m}+\alpha^m+1=\alpha^{2r}+\alpha^r+1=\alpha^r+1$, which is impossible.) 
	Let $r_1=2r-e$, then $0<r-r_1=e-r<r$ and 
 	\begin{equation*}
	0=\alpha^{2m}+\alpha^m+1=\alpha^{e(2q+1)+r_1}+\alpha^{eq+r}+1=\alpha^r+\alpha^{r_1}+1
	\end{equation*}
	On the other hand,
	\begin{equation*}
	(\alpha^{-1})^{2m}+(\alpha^{-1})^m+1=\alpha^{-1}+\alpha^{-r_1}+1=0
	\end{equation*}
 	and thus
 	\begin{equation*}
	\alpha^r+\alpha^{r-r_1}+1=0
	\end{equation*}
	From this we get $\alpha^{r-r_1}=\alpha^{r_1}$, that is, $\alpha^{|r-2r_1|}=1$ which means $e$ divides $|r-2r_1|$. 
	Since  $|r-2r_1|<e, r=2r_1$. Therefore $f$ divides the trinomial $x^{2r_1}+x^{r_1}+1$. And then we have also 
	$\alpha^{3r}=1(\alpha^{3r_1}=1)$, which implies that $e$ divides $3r(3r_1)$. Since $2r(2r_1)<e$, we get $e=3r(3r_1)$. 
	If there exists another integer $m_1$ such that
 	\begin{equation*}
	\alpha^{2m_1}+\alpha^{m_1}+1=0, ~ 2m_1<e
	\end{equation*}
	then $e=3m_1$ and therefore $m_1=r(r_1)$.  $\Box$ \\

	Now we are ready to describe the condition for divisibility of self-reciprocal trinomials by a given irreducible polynomial.
	%
	%-------------- Theorem 1 --------------
	%
	\begin{thr}
	Given an irreducible polynomial $f$ over $\F_2$, $f$ divides self-reciprocal trinomials if and only if the order of $f$ is 
	a multiple of 3.
	\end{thr}
	\textit{Proof}. Suppose $f$ divides self-reciprocal trinomials. By Lemma 1, $f$ divides self-reciprocal trinomial 
	$x^{2m}+x^m+1$ with $2m<e$ where $e$ is the order of $f$. Let $\alpha$ be any root of $f$ then 
	$\alpha^{2m}+\alpha^m+1=0$ and we get $e=3m$ as in the proof of lemma 1. Conversely suppose $e=3m$ 
	for a positive integer $m$. Let $\alpha$  be a root of $f$ then $\alpha^e=1$ that is $0=\alpha^{3m}-1=
	(\alpha^m-1)(\alpha^{2m}+\alpha^m+1)$. Since $\alpha^m \neq 1, \alpha^{2m}+\alpha^m+1=0$ and thus $f$ 
	divides the trinomial $x^{2m}+x^m+1$. $\Box$ \\

	Below we show a factorization of an arbitrary self-reciprocal trinomial over $\F_2$.
	%
	%-------------- Theorem 2 --------------
	%
	\begin{thr}
	For any odd number $m$,
 	\begin{equation*}
	x^{2m}+x^m+1=\prod_{\substack{n | m \\ 3n \nmid m}} Q_{3n}
	\end{equation*}
	where $Q_{3n}$ is the $3n$th cyclotomic polynomial over $\F_2$.
	\end{thr}
	\textit{Proof}. Suppose $n | m, 3n \nmid m$ and let $f$ be an irreducible polynomial of order $3n$ and $\alpha$ be 
	any root of $f$ in a certain extension of $\F_2$. Then $\alpha^{3n}=1$ and therefore
 	\begin{equation*}
	\alpha^{3m}-1=(\alpha^m-1)(\alpha^{2m}+\alpha^m+1)=0.
	\end{equation*}
	Since $3n \nmid m, \alpha^m-1 \neq 0$ and thus $\alpha^{2m}+\alpha^m+1=0$. Therefore $f$ divides the trinomial 
	$x^{2m}+x^m+1$. Since $Q_{3n}$ is a product of all irreducible polynomials of order $3n$, it divides the trinomial 
	$x^{2m}+x^m+1$. From $\textnormal{deg}(Q_{3n})=\phi(3n)$, it is sufficient to show
 	\begin{equation*}
	\sum_{\substack{n | m \\ 3n \nmid m}} \phi(3n)=2m
	\end{equation*}
	Using the formula $\sum_{d|n}\phi(d)=n$, we get
 	\begin{eqnarray*}
	\sum_{\substack{n | m \\ 3n \nmid m}} \phi(3n) & = & \sum_{n | m} \phi(3n) - \sum_{3n | m} \phi(3n) \\ 
	& = & \sum_{3n | 3m} \phi(3n) - \sum_{3n | m} \phi(3n) = 3m-m=2m.
	\end{eqnarray*}
	This completes the proof. $\Box$ 
	%
	%-------------- Corollary 1 --------------
	%
	\begin{cor}
	If $m$ is an odd number and $m=3^k\cdot n, ~ 3 \nmid n$ for a nonnegative integer $k$, then the self-reciprocal 
	irreducible trinomial $x^{2 \cdot 3^k}+x^{3^k}+1$ divides $x^{2m}+x^m+1$.
	\end{cor}
	\textit{Proof}. The trinomial $x^{2 \cdot 3^k}+x^{3^k}+1$ divides $Q_{3^{k+1}}$ since it is an irreducible polynomial 
	of order $3^{k+1}$. Recalling $3^k \vert m, 3^{k+1} \nmid m$, we get a desired result from Theorem 2. $\Box$\\

	We can extend Theorem 2 to any positive degree m.	
	%
	%-------------- Corollary 2 --------------
	%
	\begin{cor}
	Suppose that $m=2^k \cdot n, ~ 2 \nmid n$. Then
	\begin{equation*}
	x^{2m}+x^m+1 = \left( \prod_{\substack{n_1 | n \\ 3n_1 \nmid n}}Q_{3n_1} \right)^{2^k}.
	\end{equation*}
	\end{cor}
	\textit{Proof}. Since
	\begin{equation*}
	x^{2m}+x^m+1 = \left( x^{2n} \right)^{2^k}+\left( x^n \right)^{2^k}+1= \left( x^{2n}+x^n+1 \right)^{2^k},
	\end{equation*}
	the assertion is followed from Theorem 2. $\Box$ \\

	If an irreducible polynomial $f$ of order $e$ divides a trinomial $x^n+x^k+1$, then for all positive integer $r$ and $s, f$ 
	divides $x^{n+re}+x^{k+se}+1$ and it divides at least one trinomial of degree $<e$. Consider a number of trinomials 
	of degree $<e$ which are divided by a given irreducible polynomial. Denote as $N_f$ the number of trinomials of degree 
	$<e$ which are divided by given irreducible polynomial $f$ of order $e$.
	%
	%-------------- Theorem 3 --------------
	%
	\begin{thr}
	Let $f(x)$ be an irreducible polynomial of order $e$ which divides trinomials over $\F_2$. Then
 	\begin{equation*}
	N_f=\frac{1}{2} \textnormal{deg}\left(\textnormal{gcd}\left( 1+x^e, 1+(1+x)^e \right) \right),
	\end{equation*}
	where $\textnormal{deg}$ means the degree of the polynomial.
	\end{thr}
	\textit{Proof}. Let
	\begin{equation*}
	1+x^e=g_1(x) \cdot g_2(x) \cdot \cdots \cdot g_t(x)
	\end{equation*}
 	be a product of all irreducible polynomials whose orders divide $e$. Then we get
 	\begin{equation*}
	1+(1+x)^e=g_1(x+1) \cdot g_2(x+1) \cdot \cdots \cdot g_t(x+1).
	\end{equation*}
	Let $\alpha$ be a root of $f(x)$ then $1, \alpha, \alpha^2, \cdots, \alpha^{e-1}$ are all roots of $g_1(x), g_2(x), 
	\cdots, g_t(x)$ and $0, 1+\alpha, 1+\alpha^2, \cdots, 1+\alpha^{e-1}$ are all roots of $g_1(x+1), g_2(x+1), 
	\cdots, g_t(x+1)$. From the assumption there exists at least one pair $(i, j)$ such that $1 \leq i, j <e, i \neq j, 
	\alpha^i=\alpha^j+1$. It can be easily seen that the number of such pairs is equal to the number of common roots 
	of $1+x^e$ and $1+(1+x)^e$ that is the degree of the polynomial $\textnormal{gcd}\left( 1+x^e, 1+(1+x)^e \right)$. 
	(Note that $\textnormal{gcd}\left( 1+x^e, 1+(1+x)^e \right)$ cannot has any multiple root.) Since the different pairs 
	$(i, j)$ and $(j, i)$ correspond the same trinomial, the result is true. $\Box$ 
	%
	%-------------- Corollary 3 --------------
	%
	\begin{cor}
	The number of trinomials of degree $<2^k-1$ which are divided by a given primitive polynomial of degree $k$ is 
	$2^{k-1}-1$.
	\end{cor}

	In particular it is interesting the case when the number $N_f$ is 1.
	%
	%-------------- Theorem 4 --------------
	%
	\begin{thr}
	If $N_f$ is 1, then $f$ divides a self-reciprocal trinomial.
	\end{thr}
	\textit{Proof}. Let $e$ be an order of $f$. From Theorem 1, it is sufficient to prove that $e$ is divided by 3. 
	Suppose that $e$ is not divided by 3 and $f$ divides a trinomial $x^n+x^k+1$. Then by Theorem 1 $n \neq 2k$. 
	Let $\alpha$ be a root of $f$. Then $\alpha^{-1}$ is a root of $f^*$, the reciprocal of $f$. Since $f^*$ divides 
	$x^n+x^{n-k}+1$, 
	\begin{equation*}
	\alpha^{-n}+\alpha^{-(n-k)}+1=0,
	\end{equation*}
	that is,
	\begin{equation*}
	\alpha^{e-n}+\alpha^{e-n+k}+1=0.
	\end{equation*}
	Here
	\begin{equation*}
	0<e-n, e-n+k<e, e-n \neq e-n+k.
	\end{equation*}
	Therefore $f$ divides the trinomial $x^{e-n}+x^{e-n+k}+1$. Since $e$ is odd, $e-n \neq n$. Assume now $e-n=k$. 
	We then get
 	\begin{equation*}
	\alpha^{n+k}=\alpha^e=1.
	\end{equation*}
	Multiplying $\alpha^k$ on both sides of the equation
	\begin{equation*}
	\alpha^n+\alpha^k+1=0,
	\end{equation*}
	we have
 	\begin{equation*}
	\alpha^{2k}+\alpha^k+1=0
	\end{equation*}
	which says that $f$ divides some self-reciprocal trinomial that contradicts the assumption. Thus $f$ divides two different 
	trinomials $x^n+x^k+1$ and $x^{e-n}+x^{e-n+k}+1$ of degree $<e$, that is, $N_f \geq 2$. $\Box$

	%
	%% Section 3 Divisibility of trinomials  $x^{am}+x^{bs}+1$  %%
	%
	%
	\section{Divisibility of trinomials  $x^{am}+x^{bs}+1$}

	In this section we consider the conditions for divisibility of trinomials $x^{am}+x^{bs}+1$ by a given irreducible polynomial 
	over $\F_2$. Let $f$ be an irreducible polynomial of degree $n$ over $\F_2$ and $a$ and $b$ be positive integers. 
	In \cite{che} it was proved that if there exist positive integers $m$ and $s$ such that $f$ divides $x^{am}+x^{bs}+1$, 
	then $a$ and $b$ are not divisible by $2^n-1$. Below we give a refinement of this result.
	%
	%-------------- Theorem 5 --------------
	%
	\begin{thr}
	Let $f$ be an irreducible polynomial of order $e>1$ over $\F_2$ and $a$ and $b$ be positive integers. If there exist 
	positive integers $m$ and $s$ such that $f$ divides trinomial $x^{am}+x^{bs}+1 (am>bs)$, then $am, bs$ and $am-bs$ 
	are not divisible by $e$.
	\end{thr}
	\textit{Proof}. Let $\alpha$ be any root of $f$ in a certain extension of $\F_2$. If $am$ is divided by $e$, then 
	$\alpha^{am}=1$, so $f$ divides a polynomial $x^{am}+1$. Since $e>1, f(0) \neq 0$, and thus $f$ does not divide 
	$x^{bs}$. Therefore $f$ cannot divide the trinomial $x^{am}+x^{bs}+1$. The case where $bs$ is divided by $e$ 
	is very similar. Suppose $am-bs$ is divided by $e$. Then in the same way as above we see easily that $x^{am-bs}+1$ is 
	divided by $f$ and thus $x^{am}+x^{bs}+1=x^{bs} \left( x^{am-bs}+1 \right)+1$ is not divisible by $f$. $\Box$ \\

	If $f$ is an irreducible polynomial of order $e$ and degree $n$ over $\F_2$, then $e$ is a divisor of $2^n-1$. 
	Thus the above theorem derives directly the result in \cite{che}. And if $a=b=1$ and $f=x^2+x+1$ then the converse 
	of Theorem 5 is also true.
	%
	%-------------- Corollary 4 --------------
	%
	\begin{cor}
	The trinomial $x^n+x^k+1 (n>k)$ is divided by $x^2+x+1$ if and only if $n, k$ and $n-k$ are not divided by 3.
	\end{cor}
	\textit{Proof}. Since the order of $x^2+x+1$ is 3, the necessity is clear from above theorem. Suppose that $n, k$ and 
	$n-k$ are not divided by 3. Then we get two cases:
 	\begin{equation*}
	n \equiv 2\pmod 3, ~~ k \equiv 1\pmod 3, ~~ n-k \equiv 1\pmod 3
	\end{equation*}
	or
 	\begin{equation*}
	n \equiv 1\pmod 3, ~~ k \equiv 2\pmod 3, ~~ n-k \equiv 2\pmod 3.
	\end{equation*}
	Let $\alpha$ be any root of $x^2+x+1$ then in the first case we have
 	\begin{equation*}
	\alpha^n+\alpha^k+1=\alpha^{3n_1+2}+\alpha^{3k_1+1}+1=\alpha^2+\alpha+1=0.
	\end{equation*}
	Hence $x^2+x+1$ divides $x^n+x^k+1$. The second case is similar. $\Box$ \\

	Finally we consider the criterion for testing if an irreducible polynomial divides trinomials of type $x^{am}+x^{bs}+1$ 
	over $\F_2$.
	%
	%-------------- Theorem 6 --------------
	%
	\begin{thr}
	Let $f$ be an irreducible polynomial of order $e$ and degree $n$ over $\F_2$ and $a$ and $b$ be positive integers. 
	Then $f$ divides trinomials $x^{am}+x^{bs}+1$ if and only if $\textnormal{gcd}\big(1+x^{e_1}, 1+(1+x)^{e^2}\big)$ 
	has degree greater than 1, where
	\begin{equation*}
	e_1=\frac{e}{\textnormal{gcd}(a, e)}, \quad e_2=\frac{e}{\textnormal{gcd}(b, e)}.
	\end{equation*}
	\end{thr}
	\textit{Proof}. Let $\alpha$ be any root of $f$. Then the order of $\alpha$ in the multiplicative group $\F_{2^n}^*$ is $e$ 
	and $1, \alpha, \alpha^2, \cdots, \alpha^{e-1}$ are distinct roots of $x^e-1$. Since
	\begin{equation*}
	x^e-1=\prod_{d | n}Q_d
	\end{equation*}
 	for every $i(0 \leq i \leq e-1)$, $\alpha^i$ is a root of an irreducible polynomial whose order is a divisor of $e$. 
	In particular, $\alpha^a$ has order $e_1=\frac{e}{\textnormal{gcd}(a, e)}$ and $\alpha^a, \alpha^{2a}, 
	\cdots, \alpha^{(e_1-1)a}$ are all roots of $C_{e_1}(x):=\frac{x^{e_1}-1}{x-1}$. Similarly $\alpha^b, \alpha^{2b}, 
	\cdots, \alpha^{(e_2-1)b}$  are all roots of $C_{e_2}(x):=\frac{x^{e_2}-1}{x-1}$ and thus $1+\alpha^b, 1+\alpha^{2b}, 
	\cdots, 1+\alpha^{(e_2-1)b}$ are all roots of  $C_{e_2}(x+1)$. Hence $\alpha$ is a root of trinomial $x^{am}+x^{bs}+1$ 
	if and only if $C_{e_1}(x)$ and $C_{e_2}(x+1)$ have common root. This is equivalent to the fact that $\textnormal{gcd}
	\big(1+x^{e_1}, 1+(1+x)^{e^2}\big)$ has degree greater than 1. $\Box$ \\
	
	Put $a=b=1$ in Theorem 6. Then we have Welch's criterion.
	%
	%-------------- Corollary 5 --------------
	%
	\begin{cor} \textnormal{(\cite{gol})}
	For any odd integer $e$, the irreducible polynomials of order $e$ divide trinomials if and only if 
	$\textnormal{gcd}\big(1+x^e, 1+(1+x)^e\big)$ has degree greater than 1.
	\end{cor}


\begin{thebibliography}{20}
	
	\bibitem{bla}
	I. F. Blake, S. Gao and R. J. Lambert, Construction and distribution problems for irreducible trinomials over finite fields, 
	in Applications of Finite Fields (D. Gollmann, ed.), Oxford, Clarendon Press, 1996, 19-32.

	\bibitem{bre}
	R. Brent and P. Zimmermann, Algorithms for finding almost irreducible and almost primitive trinomials, Primes 
	and Misdemeanours : Lectures in Honour of the Sixtieth Birthday of Hugh Cowie Williams, The Fields Institute, 
	Toronto, 2004, 91-102.

	\bibitem{che}
	M. Cherif, A necessary condition of the divisibility of trinomials $x^{am}+x^{bs}+1$ by any irreducible polynomial 
	of degree $r$ over GF(2), International Journal of Algebra, \textbf{2}(13) (2008), 645-648.

	\bibitem{doc}
	C. Doche, Redundant trinomials for finite fields of characteristic 2, Proceedings of ACISP 05, LNCS \textbf{3574} 
	(2005), 122-133.

	\bibitem{gol}
	S. W. Golomb and P. F. Lee, Irreducible polynomials which divide trinomials over GF(2), IEEE Transactions on 
	Information Theory, \textbf{53} (2007), 768-774.
	
	\bibitem{lid}
	R. Lidl and H. Niederreiter, Introduction to finite fields and their applications, Cambridge University Press, 1997.

	\bibitem{swa}
	R. G. Swan, Factorization of polynomials over finite fields, Pacific Journal of Mathematics, \textbf{12} (1962), 1099-1106.

	\bibitem{tro}
	J. Tromp, L. Zhang and Y. Zhao, Small weight bases for Hamming codes, Theoretical Computer Science, 
	\textbf{181}(2) (1997), 337-345.
	
	\bibitem{yuc}
	J. L. Yucas and G. L. Mullen, Self reciprocal irreducible polynomials over finite fields, Designs, Codes and Cryptography, 
	\textbf{33} (2004), 275-281.

	\end{thebibliography}
 	\end{document}